\documentclass[10pt]{amsart}
\usepackage[cp1251]{inputenc}
\usepackage[english,russian]{babel}
\usepackage{amsmath}
\usepackage{amssymb}
\usepackage{amsfonts}

\begin{document}

\sloppy
\begin{center}
{\large\bf GENERALIZED $\delta$-DERIVATIONS}\\

\hspace*{8mm}

{\large\bf Ivan Kaygorodov}\\
e-mail: kib@math.nsc.ru

{\it 
Sobolev Inst. of Mathematics\\ 
Novosibirsk, Russia\\}
\end{center}

\

\medskip

\

\begin{center} {\bf Abstract: }\end{center}                                                                    
We defined generalized $\delta$-derivations of algebra $A$ as linear mapping $\chi$ associated with usual $\delta$-derivation $\phi$ 
by the rule  $$\chi(xy)=\delta(\chi(x)y+x\phi(y))=\delta(\phi(x)y+x\chi(y))$$ for any $x,y \in A.$
We described generalized $\delta$-derivations of 
prime alternative algebras, 
prime Lie algebras and superalgebras, 
unital algebras,
and semisimple finite-dimensional Jordan superalgebras.
In this cases we proved that generalized $\delta$-derivation is a 
generalized derivation or $\delta$-derivation.
After that we described $\delta$-superderivations of superalgebras <<KKM Double>>, arising from 
prime alternative algebras, 
prime Lie algebras and superalgebras, 
unital algebras,
and semisimple finite-dimensional Jordan superalgebras.
In the end, we constructed new examples of non-trivial $\delta$-derivations of Lie algebras.

\medskip


\medskip

\section{Введение}

\sloppy

Понятие антидифференцирования алгебры, являющееся частным случаем $\delta$-дифференцирования, 
т.е. $(-1)$-дифференцированием,
рассматривалось в работах \cite{hop2,fi}. 
В дальнейшем, в работе \cite{Fil} появляется 
определение $\delta$-дифференцирования алгебры. 
Напомним, что при
фиксированном $\delta$ из основного поля $F,$ 
под $\delta$-дифференцированием
алгебры $A$ понимают линейное отображение $\phi$, удовлетворяющее
условию 
\begin{eqnarray}\label{delta}
\phi(xy)=\delta(\phi(x)y+x\phi(y))\end{eqnarray}
 для произвольных
элементов $x,y \in A.$ 
В работе \cite{Fil} описаны
$\frac{1}{2}$-дифференци\-ро\-ва\-ния произвольной первичной
алгебры Ли $A$ ($\frac{1}{6} \in F$) с невырожденной
симметрической инвариантной билинейной формой. А именно, доказано,
что линейное отображение $\phi$: $A \rightarrow A $ является
$\frac{1}{2}$-дифференцированием тогда и только тогда, когда $\phi
\in \Gamma(A)$, где $\Gamma(A)$ --- центроид алгебры $A$. Отсюда
следует, что если $A$ --- центральная простая алгебра Ли над полем
характеристики $p \neq 2,3 $ с невырожденной симметрической
инвариантной билинейной формой, то любое
$\frac{1}{2}$-дифференцирование $\phi$ имеет вид $\phi(x)=\alpha
x,$ для некоторого $\alpha \in F$. В. Т. Филиппов доказал \cite{Fill}, что любая
первичная алгебра Ли не имеет ненулевого
$\delta$-дифференцирования, если $\delta \neq -1,0,\frac{1}{2},1$.
В работе \cite{Fill} показано, что любая первичная алгебра
Ли $A$ ($\frac{1}{6} \in \Phi$) с ненулевым антидифференцированием
является 3-мерной центральной простой алгеброй над полем частных
центра $Z_{R}(A)$ своей алгебры правых умножений $R(A)$. Также в
этой работе был построен пример нетривиального
$\frac{1}{2}$-дифференцирования для алгебры Витта $W_{1},$ т.е.
такого $\frac{1}{2}$-дифференцирования, которое не является
элементом центроида алгебры $W_{1}.$ В \cite{Filll} описаны
$\delta$-дифференцирования первичных альтернативных и нелиевых
мальцевских алгебр с некоторыми ограничениями на кольцо
операторов $\Phi$. Как оказалось, алгебры из этих классов не имеют
нетривиальных $\delta$-дифференцирований. 

В работе \cite{kay} было дано описание $\delta$-дифференцирований
простых конечномерных йордановых супералгебр над алгебраически замкнутым полем характеристики нуль.
В дальнейшем, в работе \cite{kay_lie} были описаны $\delta$-дифференцирования классических супералгебр Ли.
Работа \cite{kay_lie2} посвящена описанию $\delta$-дифференцирований 
полупростых конечномерных йордановых супералгебр над произвольным полем характеристики отличной от 2 и 
$\delta$-(супер)дифференцирований простых конечномерных
лиевых и йордановых супералгебр над алгебраически замкнутым полем характеристики нуль.
Для алгебр и супералгебр из работ \cite{kay,kay_lie,kay_lie2} было показано отсутствие
нетривиальных $\delta$-(супер)дифференцирований.
В дальнейшем, результаты \cite{kay_lie} получили обобщение в работе П. Зусмановича \cite{Zus}.
Им было дано описание $\delta$-(супер)дифференцирований первичных супералгебр Ли.
А именно, он доказал, что первичная супералгебра Ли не имеет нетривиальных $\delta$-(супер)дифференцирований  при $\delta\neq -1,0,\frac{1}{2},1.$
П. Зусманович показал, что для совершенной (т.е., такой что $[A,A]=A$) супералгебры Ли $A$ 
с нулевым центром и невырожденной суперсимметрической инвариантной билинейной
формой, у которой $A=[A,A],$   пространство $\frac{1}{2}$-(супер)дифференцирований
совпадает с (супер)центроидом супералгебры $A$. 
Также, П. Зусманович, в случае положительной характеристики поля, дал положительный ответ на вопрос
В. Т. Филип\-по\-ва о существовании делителей нуля в кольце $\frac{1}{2}$-дифференцирований первичной алгебры Ли,
сформулированный в \cite{Fill}. 
В свое время, И. Б. Кайгородовым и В. Н. Желябиным \cite{kay_zh}, рассматривались 
$\delta$-(супер)дифференцирования простых унитальных супералгебр йордановой
скобки, где ими было показано отсутствие нетривиальных $\delta$-(супер)дифференцирований
простых супералгебр йордановой скобки, не являющихся супералгебрами векторного типа и
было приведено описание $\delta$-(супер)дифференцирований
простых конечномерных унитальных йордановых супералгебр над
алгебраическим замкнутым полем характеристики $p \neq 2$.
Как следствие, была обнаружена связь между наличием нетривиальных $\delta$-дифференцирований простых унитальных супералгебр йордановых скобок 
и специальностью этой супералгебры. 
В дальнейшем, $\delta$-супердифференцирования обобщенного дубля Кантора, построенного 
на первичной ассоциативной алгебре, рассматривались в работе \cite{kay_ob_kant}. 
Цикл статей по описанию $\delta$-(супер)дифференцирований простых конечномерных йордановых супералгебр 
заканчивается работой \cite{kg_ss}, где было дано полное описание $\delta$-(супер)дифференцирований
полупростых конечномерных йордановых супералгебр над алгебраически замкнутым полем характеристики отличной от 2. 
В частности, были построены примеры нетривиальных $\frac{1}{2}$-(супер)дифференцирований для простых 
неунитальных конечномерных йордановых супералгебр. 

\
\newpage
\section{Об обобщенных $\delta$-дифференцированиях алгебр и супералгебр.} 	

\

Пусть $F$ --- поле характеристики отличной от 2. 
Напомним определение супералгебры. 
Алгебра $G$ над полем $F$ называется супералгеброй (или $\mathbb{Z}_{2}$-градуированной
алгеброй), если она представима в виде
$G=G_{0} \oplus G_{1}$, при этом справедливы
соотношения $G_{i}G_{j} \subseteq G_{i+j(mod 2)},
i,j=0,1.$ 

При фиксированном элементе $\delta \in F$, для супералгебры $A=A_0 \oplus A_1$, 
однородное линейное отображение $\phi : A \rightarrow A$ 
будем называть четным \textit{$\delta$-супердиф\-фе\-рен\-ци\-ро\-ванием}, если $\phi(A_i) \subseteq A_i$ и для однородных 
$x,y \in A_0 \cup A_1$ выполнено \begin{eqnarray*} \phi (xy)&=&\delta(\phi(x)y+x\phi(y)).\end{eqnarray*}

Под центроидом $\Gamma(A)$ супералгебры $A$ мы будем понимать множество линейных отображений 
$\chi: A \rightarrow A,$ что для произвольных элементов $a,b$ верно	
$$\chi(ab)=\chi(a)b=a\chi(b).$$

Ясно, что 1-супердифференцирование является обычным супердифференцированием; 
0-супердифференцированием является произвольный эндоморфизм $\phi$ супералгебры $A$ такой, что $\phi(A^{2})=0$. 
Ненулевое $\delta$-супердифференцирование будем считать нетривиальным, если $\delta\neq 0,1$ и $\phi \notin \Gamma(A).$
Легко видеть, что четное $\delta$-супердифференцирование будет являться $\delta$-дифференцированием. 

Пусть $A$ --- алгебра над $F$ c умножением $ab$ и
обладающая дополнительной билинейной операцией $\{ \ ,\ \}: A \times A \rightarrow A$. 
Через $Ax$ обозначим изоморфную копию алгебры $A$ и на прямой сумме векторных пространств $A\oplus Ax$
зададим умножение $\cdot$ по следующему правилу
$$a\cdot b=ab, a \cdot (bx)=(ab)x, (ax)\cdot  (bx)=\{a,b\}, \mbox{ где } a,b \in A.$$
Мы получим структуру супералгебры на $B:$ $B_0=A, B_1=Ax.$
Примерами таких супералгебр являются супералгебры, построенные по процессу Кантора (см. \cite{kay_ob_kant,kant}).
Данный процесс мы будем называть обобщенным процессом Кантора. 

При рассмотрении четных $\delta$-супердифференцирований супералгебры $B$ мы приходим к понятию обобщенного 
$\delta$-дифференцирования алгебры $A$.
Для этого достаточно заметить, что для $\phi$ --- четного $\delta$-супер\-диффе\-ренци\-рования $B$ верно 
\begin{eqnarray}\label{amod_delta}
\delta\phi(ax)b+\delta (ax)\phi (b)=\phi((ax)b)=\phi(a(bx))=\delta\phi(a)(bx)+\delta a\phi(bx).
\end{eqnarray}

В силу четности $\delta$-супердифференцирования $\phi$, мы можем положить, что $\phi(bx)=\chi(b)x,$
где $\chi \in End(A).$
Таким образом, соотношение (\ref{amod_delta}) преобразуется к выражению
\begin{eqnarray}\label{obob_delta}\chi(ab)=\delta\chi(a)b+\delta a\phi(b)=\delta\phi(a)b+\delta a\chi(b).\end{eqnarray}

Далее для алгебры $A$ линейное  отображение $\chi$, связанное с $\delta$-диффе\-ренци\-рованием $\phi$ посредством соотношения (\ref{obob_delta}),
мы будем называть {\it обобщенным $\delta$-дифференцированием}. 
Обобщенные $\delta$-дифференцирования неявно возникают в работе \cite{kay_ob_kant} при рассмотрении $\delta$-супердифференцирований 
обобщенного дубля Кантора, построенного на первичной ассоциативной алгебре.
                                                                                                  
Отметим, что обобщенное $\delta$-дифференцирование тесно связано с обобщенным дифференцированием. 
Под обобщенным дифференцированием подразумевают линейное отображение $\sigma$ алгебры $A,$ которое связано с 
некоторым дифференцированием $D$ алгебры $A$ посредством соотношения
$$\sigma(ab)=\sigma(a)b+a D(b).$$
Примером обобщенного дифференцирования, не являющегося обыкновенным дифференцированием,
 может служить отображение вида $D+\psi$, где $D$ --- дифференцирование и $\psi$ --- элемент центроида алгебры.
Заметим, что обобщенные дифференцирования рассматривались, к примеру, в работах \cite{Hvala98,Lee_shiue_01}.

Далее, во всех леммах этого  раздела, мы будем подразумевать, что $\chi$ --- обобщенное $\delta$-дифференцирование, 
связанное с $\delta$-диффе\-ренци\-рованием $\phi$, и $\chi_{\phi}=\chi-\phi$. 
Все алгебры будут рассматриваться над кольцом характеристики отличной от 2, 
а супералгебры над полем характеристики отличной от 2.

\

\textbf{Лемма 1.} Пусть $\chi$ --- обобщенное $\delta$-дифференцирование (супер)алгебры $A$, тогда $\chi_{\phi}$ --- является
$\frac{\delta}{2}$-дифференцированием $A$ и $\chi_{\phi}(ab)=\delta a\chi_{\phi}(b)=\delta\chi_{\phi}(a)b.$

\

\textbf{Доказательство}. Рассматривая разность между выражениями (\ref{obob_delta}) и (\ref{delta}), мы получим
$$\chi_{\phi}(ab)=\delta a\chi_{\phi}(b)=\delta\chi_{\phi}(a)b.$$
Далее, воспользовавшись полученным равенством, 
легко имеем $$\chi_{\phi}(ab)=\frac{1}{2}(\delta a\chi_{\phi}(b)+\delta\chi_{\phi}(a)b).$$
Данное означает, что $\chi_{\phi}$ является $\frac{\delta}{2}$-дифференцированием алгебры $A.$ Лемма доказана.

\

Ясно, что обобщенное 1-дифференцирование является отображением вида $D+\psi,$ где $D$ --- дифференцирование, а $\psi$ --- элемент центроида.
Обобщенным 0-дифференцированием является произвольный эндоморфизм $\chi$ с условием $\chi(A^2)=0.$
Обобщенное $\delta$-дифференцирование $\chi$ является нетривиальным, если $\delta \neq 0,1$ 
и $\chi$ не является $\delta$-дифференцированием.
Следует отметить, что условие $\chi_{\phi}=0$, непосредственно, влечет тривиальность $\chi$.

\

Напомним, что первичной (супер)алгеброй называют алгебру $A$, 
которая не обладает двумя взаимно аннулирующими друг друга левыми идеалами.
В частности, для первичной (супер)алгебры $A$ верно $$Ann(A)=\{ x\mbox{ | }Ax=0 \} =\{0\}.$$

\

\textbf{Теорема 2.} Первичная алгебра Ли $A$ не имеет нетривиальных обобщенных $\delta$-дифференцирований. 

\

\textbf{Доказательство.} В силу показанного в лемме 1, мы можем заключить, что выполняется следующая цепочка соотношений

$$\delta(x\chi_{\phi}(z))y+\delta x(y\chi_{\phi}(z))=\delta(xy)\chi_{\phi}(z)=\chi_{\phi}((xy)z)=$$
$$\chi_{\phi}((xz)y+x(yz))=\delta^2((x\chi_{\phi}(z))y+x(y\chi_{\phi}(z))).$$
Откуда имеем 
$$0=(x\chi_{\phi}(z))y+x(y\chi_{\phi}(z))=\chi_{\phi}(xy)z.$$
Таким образом, 
$$\chi_{\phi}(A^2) \subseteq Ann(A)= \{ 0\}.$$
Отсюда получаем $\chi_{\phi}(x)y=0$ и 
$$\chi_{\phi}(A) \subseteq Ann(A)= \{ 0\}.$$
Что эквивалентно тривиальности $\chi.$ 
Теорема доказана.

\

\textbf{Теорема 3.} Первичная лиева супералгебра $A$ не имеет нетривиальных обобщенных $\delta$-дифференцирований. 

\

\textbf{Доказательство.} 
Легко понять, что пространство $End(A)$ является $\mathbb{Z}_2$-градуированным, 
то есть, любое линейное отображение $\psi \in End(A)$ мы можем представить в виде суммы четного и нечетного отображений $\psi_0+\psi_1$, где 
$\psi_0(A_i) \subseteq A_i$ и $\psi_1(A_i) \subseteq A_{i+1}.$

Будем считать, что $\chi$ и $\phi$ являются четными отображениями, то есть  верно 
$$\chi(A_i)\subseteq A_i, \phi(A_i)\subseteq A_i.$$
Тогда, в силу показанного в лемме 1, мы можем заключить, что для однородных элементов $x,y,z \in A_0 \cup A_1,$ 
выполняется следующая цепочка соотношений
$$\delta(x\chi_{\phi}(z))y+(-1)^{p(x)p(z)}\delta x(y\chi_{\phi}(z))=\delta(xy)\chi_{\phi}(z)=\chi_{\phi}((xy)z)=$$
$$\chi_{\phi}((xz)y+(-1)^{p(x)p(z)}x(yz))=\delta^2((x\chi_{\phi}(z))y+(-1)^{p(x)p(z)}x(y\chi_{\phi}(z))).$$
Откуда получаем
$$0=(x\chi_{\phi}(z))y+(-1)^{p(x)p(z)}x(y\chi_{\phi}(z))=(xy)\chi_{\phi}(z).$$
Понятно, что из законов дистрибутивности, также следует, 
что равенство $\chi_{\phi}(xy)z=0$ выполняется для произвольных $x,y,z \in A,$
где $\chi_{\phi}$ --- четное отображение, определенное выше.

Пусть $\chi$ и $\phi$ являются нечетными отображениями, то есть 
верно $$\chi(A_i)\subseteq A_{i+1}, \phi(A_i)\subseteq A_{i+1}.$$
Положим $x_i,y_i \in A_i$ и $x,y,z \in A.$

Легко заметить, что 
$$\delta(x_0\chi_{\phi}(z))y+\delta x_0(y\chi_{\phi}(z))=\delta(x_0y)\chi_{\phi}(z)=\chi_{\phi}((x_0y)z)=$$
$$\chi_{\phi}((x_0z)y+x_0(yz))=\delta^2((x_0\chi_{\phi}(z))y+x_0(y\chi_{\phi}(z))).$$
Откуда имеем 
$$0=(x_0\chi_{\phi}(z))y+x_0(y\chi_{\phi}(z))=(x_0y)\chi_{\phi}(z),$$
то есть $\chi_{\phi}(x_0y)z=0.$

Отметим, что 
$$\chi_{\phi}(x_1y_0)=\delta\chi_{\phi}(x_1)y_0=-\delta y_0 \chi_{\phi}(x_1)=-\delta \chi_{\phi}(y_0)x_1=-\delta x_1 \chi_{\phi}(y_0)=-\chi_{\phi}(x_1y_0),$$
то есть $\chi_{\phi}(x_1y_0)z=0.$

Заметим, что 
$$\chi_{\phi}(x_1y_1)=\delta\chi_{\phi}(x_1)y_1=-\delta y_1 \chi_{\phi}(x_1)=-\chi_{\phi}(y_1x_1)=-\chi_{\phi}(x_1y_1),$$
то есть  $\chi_{\phi}(x_1y_1)z=0$ и  $\chi_{\phi}(x_1y)z=0.$

Теперь мы можем заключить, что  $\chi_{\phi}(xy)z=0,$
где $x,y,z$ произвольные элементы $A$ и $\chi_{\phi}$ либо четное, либо нечетное отображение.

Таким образом, 
$$\chi_{\phi}(A^2) \subseteq Ann(A)= \{ 0\}.$$
Откуда получаем $\chi_{\phi}(x)y=0$ и 
$$\chi_{\phi}(A) \subseteq Ann(A)= \{ 0\}.$$
Что эквивалентно тривиальности $\chi.$ 
Теорема доказана.

\

\textbf{Теорема 4.} Первичная альтернативная алгебра $A$ не имеет нетривиальных обобщенных $\delta$-дифференцирований. 

\

\textbf{Доказательство.} В силу показанного в лемме 1, мы можем заключить, что выполняется следующая цепочка соотношений
$$\delta \chi_{\phi}(y)x^2=\chi_{\phi}(yx^2)=\chi_{\phi}((yx)x)=\delta^2((\chi_{\phi}(y)x)x)=\delta^2\chi_{\phi}(y)x^2.$$
Таким образом, мы получаем $y\chi_{\phi}(x^2)=0,$ то есть $\chi_{\phi}(x^2) \subseteq Ann(A)=\{ 0\}.$
Откуда, посредством линеаризации, легко вытекает, что $\chi_{\phi}(xy+yx)=0.$
Рассмотрим на $A$ новое умножение: 
$$a\odot b=\frac{1}{2}(ab+ba).$$
Полученную алгебру с новым умножением $\odot$, как обычно, будем обозначать $A^{(+)}.$
В силу работ \cite{pch04,pch07}, если $A$ --- альтернативная первичная алгебра, то $A^{(+)}$ --- йорданова первичная алгебра.
Таким образом, мы имеем, что $$\chi_{\phi}(x) \odot y=0,$$
что влечет $$\chi_{\phi}(x) \subseteq Ann(A^{(+)}) =\{0\}.$$
Исходя из полученного, мы имеем тривиальность $\chi.$ Теорема доказана. 

\

\textbf{Теорема 5.} Унитальная (супер)алгебра $A$ не имеет нетривиальных обобщенных $\delta$-дифференцирований. 

\

\textbf{Доказательство.} 
Заметим, что если $e$ --- единица (супер)алгебры $A$, то 
$$\chi_{\phi}(e)=\chi_{\phi}(ee)=\delta\chi_{\phi}(e) \mbox{ и  }\chi_{\phi}(e)=0.$$
Таким образом, легко видеть, что $$\chi_{\phi}(x)=\chi_{\phi}(ex)=2\chi_{\phi}(e)x=0.$$
Откуда мы получаем тривиальность $\chi.$ Теорема доказана.

\

В частности, теорема 5 дает отсутствие нетривиальных обобщенных $\delta$-дифференцирований для 
полупростых конечномерных йордановых алгебр и всего класса структуризуемых алгебр над произвольным полем характеристики отличной от 2.

\

\textbf{Теорема 6.} Полупростая конечномерная йорданова супералгебра $A$ над алгебраически замкнутым полем характеристики отличной от 2
не имеет нетривиальных обобщенных $\delta$-дифференцирований. 

\

{\bf Доказательство.} 
Согласно работе \cite{Zel_ssjs},  если $A$ полупростая конечномерная йорданова супералгебра над алгебраически замкнутым полем 
характеристики отличной от 2, то 
$A=\bigoplus_{i=1}^{s}T_i\oplus J_1 \oplus \ldots \oplus J_t,$ где $J_1,\ldots,J_t$ --- простые йордановы супералгебры 
и $T_i=J_{i1}\oplus\ldots \oplus J_{ir_i}+K_i \cdot 1,$
 $K_1, \ldots, K_s$ --- расширения поля $F$ и $J_{i1},\ldots , J_{ir_i}$ --- простые 
неунитальные йордановы супералгебры над полем $K_i.$
Пусть некоторая супералгебра $J$ представима в виде прямой суммы супералгебр $B \oplus C$ 
и $b$ --- элемент супералгебры $B$, который не являются делителем нуля,
 а $c$ --- произвольный элемент $C$.
Тогда 
$$0=\chi_{\phi}(b c)=\delta b \chi_{\phi}(c),$$
откуда получаем, что $\chi_{\phi}(c) \in C.$ 
Благодаря чему, можем заключить, что
 $\chi_{\phi}(J_i) \subseteq J_i$ и $\chi_{\phi}(T_i) \subseteq T_i.$ 
Также, легко заметить, что $\phi(J_i) \subseteq J_i$ и $\phi(T_i) \subseteq T_i.$ 
Исходя из теоремы 5, мы заключаем, что ограничение $\chi$ на унитальные супералгебры $T_i$ и $J_i$ --- тривиальны.

Следовательно, нам достаточно показать тривиальность ограничения $\chi$ на $J_i$ в случае, когда $J_i$ является 
неунитальной простой конечномерной йордановой супералгеброй. 
Данные супералгебры исчерпываются супералгебрами $K_3,V_{1/2}(Z,D)$ и супералгеброй $K_9$ в случае характеристики поля $p=3$.
Их определения можно найти, к примеру, в работах \cite{kay,kg_ss,Zel_ssjs}. 
В частности, известно, что четные части супералгебр $K_3,K_9,V_{1/2}(Z,D)$ являются унитальными алгебрами.
Согласно \cite[теорема 10]{kg_ss}, в этом случае, $J_i$ не имеет нетривиальных $\delta$-дифференцирований при $\delta \neq \frac{1}{2}$.
Таким образом, по лемме 1 и определению нетривиального обобщенного $\delta$-дифференцирования, $\chi_{\phi}=0$ при $\delta\neq 2.$
Случай $\delta=2$ рассмотрим подробней. Легко понять, что если $e$ --- единица $(J_i)_0,$ то $\chi_{\phi}(e) \in (J_i)_1.$
Из определения супералгебр $K_3,V_{1/2}(Z,D)$, известно, что $2ez=z,$ при $z \in (J_i)_1.$
Следовательно, если $J_i=K_3,V_{1/2}(Z,D)$, то для $z \in (J_i)_1$ верно $\chi_{\phi}(z)=2\chi_{\phi}(ez)=4e\chi_{\phi}(z),$ 
то есть $\chi_{\phi}(z)=0$ (при $p \neq 3$) и $\chi_{\phi}(z) \in (J_i)_0$ (при $p=3$).
Если $J_i=K_3$, то известно, что для нее существуют такие $w,t \in (J_i)_1,$ что $e=wt,$ откуда, если $p\neq 3$, вытекает $\chi_{\phi}(e)=0.$
Пусть теперь $p\neq 3, J_i=V_{1/2}(Z,D),$ тогда 
$$0=\chi_{\phi}(ax)=4(\chi_{\phi}(a)\cdot x) \mbox{ и }\chi_{\phi}(a)=a_{\chi_{\phi}}x,$$
то есть $D(a_{\chi_{\phi}})=0,$ что дает $a_{\chi_{\phi}}=\alpha^{a_{\chi_{\phi}}} e, 
\alpha^{a_{\chi_{\phi}}}$ --- элемент основного поля супералгебры $J_i$.
Заметим, что $$\alpha^{a_{\chi_{\phi}}}ex=\chi_{\phi}(a)=2\chi_{\phi}(e)a=\alpha^{e_{\chi_{\phi}}}ax$$ 
и, учитывая, что $e$ и $a$ мы можем взять линейно независимые, получаем $\chi_{\phi}=0.$

Если $p=3$ и $J_i=V_{1/2}(Z,D)$ то при $z,y \in (J_i)_0,$ мы получаем $\chi_{\phi}(e)=bx$ и 
$$\chi_{\phi}(z)= 2\chi_{\phi}(e) \cdot z=(zb)x,$$
$$\chi_{\phi}(zx)=\chi_{\phi}(e) \cdot zx=D(b)z-bD(z).$$
Таким образом, имеем что 
$$\chi_{\phi}(zx \cdot yx)=2 \chi_{\phi}(zx) \cdot yx,$$
$$((D(z)y-zD(y))b)x=((D(b)z-bD(z))y)x,$$
то есть $D(zyb)=0,$ что влечет $D(b)=0$ и $b=\beta e$, где $\beta$ --- элемент основного поля супералгебры $J_i.$
Следовательно, $D(\beta z)=0$.
Замечая, что $z$ может быть линейно независимо c $e$, и
в силу того, что $D$ обнуляет только элементы вида $\gamma e$, 
где $\gamma$ --- элемент основного поля супералгебры $J_i,$
получаем $\beta=0.$
Полученное дает тривиальность $\chi.$

Если $p=3$ и $J_i=K_9,K_3$, то при $z,t \in (J_i)_1$ имеем 
$$\chi_{\phi}(zt)=2\chi_{\phi}(z)\cdot t= 2t \cdot \chi_{\phi}(z)=\chi_{\phi}(tz)=-\chi_{\phi}(zt).$$
Отметим, что здесь мы воспользовались тем, что $\chi_{\phi}(z) \in (J_i)_0.$
Таким образом, $\chi_{\phi}((J_i)_1^2)=0,$ что влечет $\chi_{\phi}(e)=0$ и тривиальность $\chi.$

Исходя из вышеприведенных рассуждений, теорема доказана.

\section{О $\delta$-дифференцированиях обобщенного дубля Кантора.}

Напомним, что в работе \cite{kay_ob_kant} было введено и рассматривалось понятие обобщенного дубля Кантора. 
Как следует из результатов работы \cite{kay_zh}, особый интерес представляет рассмотрение $\delta$-дифференцирований
обобщенного дубля Кантора, построенного на супералгебре с тривиальной нечетной компонентой. Данное заключение вытекает из того факта, что
простые унитальные супералгебры йордановой скобки имеют нетривиальные $\delta$-дифференцирования только когда они построены на
супералгебрах с тривиальной нечетной частью \cite{kay_zh}.

Пусть $A$ --- алгебра c операцией $\{ \ , \ \} : A \times A \rightarrow A$ 
и $K(A)$ --- супералгебра, полученная с помощью обобщенного процесса Кантора, описанного в параграфе 2.
Через $\Delta_{\delta}(A)$ и $\Gamma(A)$ обозначим, соответственно, множество $\delta$-дифференцирований $A$ и центроид алгебры $A$, 
а через $\Delta_{\delta}(A, \{ \ , \ \})$ --- множество $\delta$-дифференцирований $A$ по операции $\{ \ , \ \}.$
Для $\phi$ --- линейного отображения супералгебры $K(A)$, под $\phi|_{A}$ будем подразумевать ограничение отображения $\phi$ на подалгебру $A$.

\medskip
\textbf{Теорема 7.} 
Пусть $\phi$ --- нетривиальное четное $\delta$-дифференцирование супералгебры $K(A),$
где $A$ --- первичная альтернативная алгебра над полем характеристики отличной от 2,3. 
Тогда $$\delta=\frac{1}{2}\mbox{ и }
\phi(ax)= \phi|_{A}(a) x, \mbox{ где } 
\phi|_{A} \in \Gamma(A) \cap \Delta_{\frac{1}{2}}(A, \{ \ , \ \}).$$
\medskip

{\bf Доказательство.} 
Пользуясь предварительными рассуждениями из параграфа 2 и теоремой 4, 
легко понять, что при $\phi$ --- четном $\delta$-дифференцировании супералгебры $K(A)$,
мы имеем $\phi(ax)=\phi(a)x$ для любого $a \in A.$
C другой стороны, мы видим, что 
$$\phi(ax \cdot bx )= \delta \phi(ax) \cdot bx +\delta ax \cdot \phi(bx),$$
то есть 
$$\phi \{ a,b\}=\delta \{ \phi(a),b \} +\delta \{a, \phi(b)\}.$$

Таким образом, четные $\delta$-дифференцирования супералгебры $K(A)$ определяются множеством 
отображений $\Delta_{\delta}(A) \cap \Delta_{\delta}(A, \{ \ , \ \})$,
где каждое отображение $\psi \in \Delta_{\delta}(A) \cap \Delta_{\delta}(A, \{ \ , \ \})$ 
продолжается на нечетную компоненту супералгебры $K(A)$ по принципу $\psi(ax)=\psi(a)x.$

Напомним, что в силу \cite{Filll}, первичные альтернативные алгебры (над полем характеристики отличной от 2,3) 
не имеют нетривиальных $\delta$-дифференцирований, то есть 
супералгебра $K(A)$, построенная на первичной альтернативной алгебре $A$ не имеет нетривиальных четных 
$\delta$-дифференцирований при $\delta \neq \frac{1}{2}$.

Исходя из выше сказанного, множество четных $\frac{1}{2}$-дифференцирований супералгебры $K(A)$ 
определяются посредством множества $\Gamma(A) \cap \Delta_{\frac{1}{2}}(A, \{ \ ,\ \})$ 
и условия $\phi(ax)=\phi|_{A}(a)x,$ где $\phi|_{A} \in \Gamma(A) \cap \Delta_{\frac{1}{2}}(A, \{ \ ,\ \}).$
Теорема доказана.

\

\textbf{Теорема 8.} 
Пусть $\phi$ --- нетривиальное четное $\delta$-дифференцирование супералгебры $K(A),$
где $A$ --- первичная лиева алгебра (либо первичная супералгебра Ли, рассматриваемая как алгебра). 
Тогда $\delta=\frac{1}{2},-1$;
множество четных $(-1)$-дифференцирований (соответственно, четных $\frac{1}{2}$-дифференцирований) супералгебры $K(A)$ 
определяется посредством множества отображений  
$\Delta_{-1}(A) \cap \Delta_{-1}(A, \{ \ , \ \})$ (соответственно, $\Delta_{\frac{1}{2}}(A) \cap \Delta_{\frac{1}{2}}(A, \{ \ ,\ \})$) 
и условия $\phi(ax)=\phi(a)x.$

\

{\bf Доказательство.} 
Пользуясь предварительными рассуждениями из параграфа 2 и теоремой 2 (в случае супералгебры, теоремой 3), 
легко понять, что при $\phi$ --- четном $\delta$-дифференцировании супералгебры $K(A)$,
мы имеем $\phi(ax)=\phi(a)x$ для любого $a \in A.$
C другой стороны, мы видим, что 
$$\phi(ax \cdot bx )= \delta \phi(ax) \cdot bx +\delta ax \cdot \phi(bx),$$
то есть 
$$\phi \{ a,b\}=\delta \{ \phi(a),b \} +\delta \{a, \phi(b)\}.$$

Таким образом, четные $\delta$-дифференцирования супералгебры $K(A)$ определяются множеством 
отображений $\Delta_{\delta}(A) \cap \Delta_{\delta}(A, \{ \ , \ \})$,
где каждое отображение $\psi \in \Delta_{\delta}(A) \cap \Delta_{\delta}(A, \{ \ , \ \})$ 
продолжается на нечетную компоненту супералгебры $K(A)$ по принципу $\psi(ax)=\psi(a)x.$

Согласно результатам \cite{hop2,Fil,Fill} (в случае супералгебры, \cite{Zus}), 
первичные лиевые алгебры (соответственно, первичные лиевы супералгебры, рассматриваемые как алгебры) 
не имеют нетривиальных $\delta$-дифференцирований 
при $\delta \neq -1, \frac{1}{2}$. Таким образом, 
множество четных $\delta$-дифференцирований супералгебры $K(A)$, построенной на первичной лиевой алгебре $A$ 
(либо первичной лиевой супералгебре $A$, рассматриваемой как алгебра),
исчерпывается случаями $(-1)$-дифференцирований и $\frac{1}{2}$-дифференцирований.
Множество $(-1)$-дифференцирований (соответственно, $\frac{1}{2}$-дифференцирований) супералгебры $K(A)$ 
определяется посредством множества отображений  
$\Delta_{-1}(A) \cap \Delta_{-1}(A, \{ \ , \ \})$ (соответственно, $\Delta_{\frac{1}{2}}(A) \cap \Delta_{\frac{1}{2}}(A, \{ \ ,\ \})$) 
и условия $\phi(ax)=\phi(a)x.$
Теорема доказана.

\medskip 

\textbf{Теорема 9.} 
Пусть $\phi$ --- нетривиальное четное $\delta$-дифференцирование супералгебры $K(A),$
где $A$ --- унитальная алгебра или полупростая йорданова супералгебра над алгебраически замкнутым полем характеристики отличной от 2. 
Тогда $$\delta=\frac{1}{2}\mbox{ и }
\phi(ax)= \phi|_{A}(a) x, \mbox{ где } 
\phi|_{A} \in \Delta_{\frac{1}{2}}(A) \cap \Delta_{\frac{1}{2}}(A, \{ \ , \ \}).$$
{\bf Доказательство.} Утвеждение вытекает из результатов об описании $\delta$-дифференцирований 
полупростой йордановой алгебры над алгебраически замкнутым полем характеристики отличной от 2 \cite{kg_ss} и
$\delta$-дифференцирований унитальной алгебры \cite[Теорема 2.1]{kay}. 
Согласно этим результатам, на данных классах алгебр и супералгебр возможны только нетривиальные $\frac{1}{2}$-дифференцирования.
Таким образом, пользуясь схемой доказательств теорем 7 и 8, мы получаем требуемое.
Теорема доказана.

\section{Новые примеры нетривиальных $\delta$-дифференцирований алгебр Ли.}

Рассмотрим супералгебру $K(A)$, построенную по обобщенному процессу Кантора (см. параграф 2), 
где $A$ --- простая алгебра Ли с умножением $[\ ,\ ].$ 
Умножение $\cdot$ супералгебры $K(A)$ будет задаваться следующим правилом:
$$a\cdot b=[a,b],\ a\cdot bx=[a,b]x,\ ax\cdot b=[a,b]x,\ ax\cdot bx =[a,b].$$

Легко заметить, что полученная алгебра $K(A)=A \oplus Ax$ будет алгеброй Ли. 
В данном случае, множество четных $(-1)$-дифференцирований (соответственно, $\frac{1}{2}$-дифференцирований) супералгебры $K(A)$ 
определяется множеством $\Delta_{-1}(A)$ (соответственно, $\Delta_{\frac{1}{2}}(A)$) и условием $\phi(ax)=\phi(a)x$ 
для $\phi \in \Delta_{\delta}(A), \delta=-1,\frac{1}{2}.$

Исходя из вышеизложеного, теоремы 8, и известных примеров нетривиальных 
$(-1)$-дифференцирований для простой алгебры $sl_2$ (см. \cite{hop2} и пример 1 ниже) 
и нетривиальных $\frac{1}{2}$-дифференцирований для простой алгебры Витта $W_1$ (см. \cite{Fill} и пример 2 ниже),
мы получаем новые примеры нетривиальных $(-1)$-дифференцирований и $\frac{1}{2}$-дифференцирований для алгебр Ли,
построенных из простых лиевых алгебр по обобщенному процессу Кантора.

\

Напомним вышеупомянутые примеры нетривиальных $\delta$-дифферен\-ци\-ро\-ва\-ний для простых алгебр Ли. 

\

\textbf{Пример 1 \cite{hop2}.} 
Пусть умножение алгебры $sl_2 \cong k^3$ определено следующим образом 
   $$\left[\begin{array}{c}
    a  \\
    b  \\
    c  \\
    \end{array} \right] \left[\begin{array}{c}
    x  \\
    y  \\
    z  \\
    \end{array} \right] =
    \left[\begin{array}{c}
    bx - cy  \\
    2ay - 2bx  \\
    2cx - 2az  \\
    \end{array} \right], \mbox{ где }
    \left[\begin{array}{cc}
    a  & b\\
    c & -a  \\
    \end{array} \right] \mapsto \left[\begin{array}{c}
    a  \\
    b  \\
    c  \\
    \end{array} \right].$$ 
Через $Antider(sl_2)$ означим пространство антидифференцирований алгебры $sl_2$, тогда 
$Antider(sl_2)= \left\{ \left[
\begin{array}{ccc}
    -2a & b & c  \\
    2c & a & d  \\
    2b & e & a  \\
    \end{array}\right]  | a,b,c,d,e \in k     \right\}.$

\

\textbf{Пример 2 \cite{Fill}.} Алгебра $W_1$ представлет собой алгебру дифференцирований кольца многочленов от одной переменной $x$.
Отображение 
$$R: y \rightarrow ys(x_1,x_2,x_3),$$ где 
$$ys(x_1,x_2,x_3)=\sum\limits_{\sigma \in S_3} (-1)^{\sigma}yx_{\sigma(1)}x_{\sigma(2)}x_{\sigma(3)}$$
--- стандартный лиев многочлен четвертой степени, при некоторых $x_1,x_2,x_3$ --- элементах алгебры $W_1$,
будет являться нетривиальным $\frac{1}{2}$-дифференцированием.

\medskip

В заключение, автор выражает благодарность В. Н. Желябину и А. П. Пожидаеву за внимание к работе и конструктивные замечания.


\newpage

\begin{thebibliography}{2}

\bibitem{hop2} Hopkins N. C., {\it Generalizes Derivations of Nonassociative
Algebras}, Nova J. Math. Game Theory Algebra, \textbf{5} 
(1996), \No 3, 215--224.

\bibitem{fi} Филиппов В. Т., {\it Об алгебрах Ли, удовлетворяющих тождеству 5-ой степени}, Алгебра и логика, {\bf 34} (1995), \No 6, 681--705.

\bibitem{Fil} Филиппов~В.~Т.,
\textit{О $\delta$-дифференцированиях алгебр Ли}, Сиб. матем. ж.,
\textbf{39} (1998), \No 6, 1409--1422.

\bibitem{Fill} Филиппов~В.~Т.,
\textit{О $\delta$-дифференцированиях первичных алгебр Ли}, Сиб. матем. ж.,
\textbf{40} (1999), \No 1, 201--213.

\bibitem{Filll} Филиппов~В.~Т.,
\textit{О $\delta$-дифференцированиях первичных альтернативных и
мальцевских алгебр}, Алгебра и Логика, \textbf{39} (2000), \No 5,
618--625.  


\bibitem{kay} Кайгородов~И.~Б., 
\textit{О $\delta$-дифференцированиях простых конечномерных йордановых супералгебр}, 
Алгебра и логика \textbf{46} (2007), \No 5, 585--605. [ http://arxiv.org/abs/1010.2419 ]

\bibitem{kay_lie} Кайгородов~И.~Б., 
\textit{О $\delta$-дифференцированиях классических супералгебр Ли}, Сиб. матем. ж. \textbf{50} (2009), \No 3, 547--565. [ http://arxiv.org/abs/1010.2807 ]

\bibitem{kay_lie2} Кайгородов~И.~Б., 
\textit{О $\delta$-супердифференцированиях простых конечномерных йордановых и лиевых супералгебр}, Алгебра и логика \textbf{49} (2010), \No 2, 195--215. [ http://arxiv.org/abs/1010.2423 ]

\bibitem{kay_ob_kant} Кайгородов И. Б., \textit{Об обобщенном дубле Кантора}, Вестник Самарского гос. университета {\bf 78}  (2010), \No 4, 42--50. [ http://arxiv.org/abs/1101.5212 ]

\bibitem{Zus} Zusmanovich~P., 
\textit{On $\delta$-derivations of Lie algebras and superalgebras}, J. of Algebra {\bf 324} (2010), \No 12, 3470--3486. 
[ http://arxiv.org/abs/0907.2034 ]

\bibitem{kay_zh} Желябин В. Н., Кайгородов И. Б., 
\textit{О $\delta$-супердифференцированиях простых супералгебр йордановой скобки}, 
Алгебра и анализ, {\bf 23} (2011), \No 4, 40--58. [ http://arxiv.org/abs/1106.2884 ]

\bibitem{kg_ss} Кайгородов И. Б., 
{\it О $\delta$-супердифференцированиях полупростых конечномерных йордановых супералгебр}, Мат. заметки, принято к печати, [ http://arxiv.org/abs/1106.2680 ]


\bibitem{kant} Кантор~И.~Л., 
\textit{Йордановы и лиевы супералгебры, определяемые алгеброй Пуассона}, в сб. <<Алгебра и анализ>>, Томск, изд-во
ТГУ (1990), 89--126.

\bibitem{Hvala98} Hvala B., {\it Generalized derivations in rings,}
Comm. Algebra, {\bf 26} (1998), \No 4, 1147--1166. 

\bibitem{Lee_shiue_01} Lee T.-K., Shiue W-K., {\it Identities with generalized derivations,}
Comm. Algebra, {\bf 29} (2001), \No 10, 4437--4450. 

\bibitem{pch04} Пчелинцев С. В., {\it Первичные альтернативные алгебры, близкие к коммутативным}, 
Изв. РАН. Сер. матем., {\bf 68} (2004), \No 1, 183–206. 

\bibitem{pch07} Пчелинцев С. В., {\it Исключительные первичные альтернативные алгебры}, 
Сиб. матем. журн., {\bf 48} (2007), \No 6, 1322–1337. 

\bibitem{Zel_ssjs} Zelmanov E., \textit{Semisimple finite dimensional Jordan superalgebras}, (English)
[A] Fong, Yuen (ed.) et al., Lie algebras, rings and related topics. Papers of the 2nd Tainan-Moscow 
international algebra workshop '97, Tainan, Taiwan, January 11--17, 1997. Hong Kong: Springer (2000), 227--243. 




\end{thebibliography}
\end{document}